\documentclass[12pt]{article}
\usepackage{graphicx}
\usepackage{amsfonts}
\usepackage{amssymb}
\usepackage{amsmath}
\usepackage{amscd}

\usepackage{amsthm}

\usepackage[colorlinks,urlcolor=blue]{hyperref}

\newtheorem*{remark}{Remark}
\newtheorem{theorem}{Theorem}[section]
\newtheorem{lemma}[theorem]{Lemma}
\newtheorem{proposition}[theorem]{Proposition}

\begin{document}

\title{\bf Application of Lagrangian mechanics equations for finding of the minimum distance between smooth submanifolds in N-dimensional Euclidean space -- Part II}
\providecommand{\keywords}[1]{\textbf{\textit{Keywords: }} #1}

\author{
  Stanislav S. Zub\\
  Faculty of Computer Science and Cybernetics,\\
  Taras Shevchenko National University of Kyiv,\\
  Glushkov boul., 2, corps 6.,\\
  Kyiv, Ukraine 03680\\
  \texttt{stah@univ.kiev.ua}\\
  \\
  Sergiy I. Zub\\
  Institute of Metrology,\\
  Mironositskaya st., 42,\\
  Kharkiv, Ukraine 61002\\
  \texttt{sergii.zub@gmail.com}\\
  \\
  Vladimir V. Semenov\\
  Faculty of Computer Science and Cybernetics,\\
  Taras Shevchenko National University of Kyiv,\\
  Glushkov boul., 2, corps 6.,\\
  Kyiv, Ukraine 03680\\
  \texttt{semenov.volodya@gmail.com}}

\maketitle
\thispagestyle{empty}

\newpage
\begin{abstract}
The method of finding the minimal distance between smooth non crossing submanifolds in N-dimensional Euclidean space are presented. It based on solution of the equations that describe the dynamics of the pair of material points.
The dynamical system can be presented as a natural mechanical system determined by Riemannian geometry on the manifold and chosen potential energy. Such an approach makes it possible to find Lyapunov function of the considered system and to formulates the requirements on the form of potential energy that brings to the convergence of the method.

.

\end{abstract}

\newpage
\tableofcontents

%
\newcommand{\s}[1]{\ensuremath{\boldsymbol{\sigma}_{#1}}}
\newcommand{\bsym}[1]{\ensuremath{\boldsymbol{#1}}}
\newcommand{\w}{\ensuremath{\boldsymbol{\wedge}}}
\newcommand{\lc}{\ensuremath{\boldsymbol{\rfloor}}}
\newcommand{\rc}{\ensuremath{\boldsymbol{\lfloor}}}

\thispagestyle{empty}

\newpage

\bigskip
\section{ Introduction }

\bigskip
Application of Lagrange mechanics equations (see section \ref{secLagrangeForm}) for finding out the minimal distance between smooth non crossing  surfaces\footnote[1]{{\it surface} is used to denote an submanifold of any sub dimension in N-dimensional Euclidean space} in N-dimensions Euclidean space are presented.

Method base on solving of equations that describe the dynamics of the material point pair that moving under the action of a mutual force of attraction between them. Each point on the own surface is held back by holonomic constraints. The potential energy of points interaction depends only from the distance (see article \cite{GOPM}).
For stopping-down in position that corresponds to minimum distance we add dissipation through Rayleigh function.
The correspondent equations follow from d'Alembert principle and represented in section \ref{secDalamber}.

Because of every of the surfaces are Riemannian manifold, so their direct product is also a Riemannian manifold
(see section \ref{RimanGeom}) and our dynamical system can be represented as a natural mechanical system \cite{SteveSmale67}
that is defined of Riemannian geometry on manifold and potential energy \cite[(1.1.7),p.3]{MarRat98}.

Such an approach makes it possible to find out Lyapunov function of considered system and formulates the requirements on the form of potential energy that gives the method convergence (see section \ref{Lyapunov}).

\bigskip
\section{ Set up of the problem and notations }

Let $\mathfrak{N}$, $\mathfrak{M}$ are smooth non crossing surfaces in Euclidean space $E^N$.

In local coordinate systems the points $\xi$,$\eta$ of surfaces are described by radius-vector
\begin{equation}
\label{Rvectors}
   \begin{cases}
      \vec{x}(\xi^1,\cdots,\xi^n)\in\mathfrak{N},\quad n<N;\\
      \vec{y}(\eta^1,\cdots,\eta^m)\in\mathfrak{M},\quad m<N.
   \end{cases}
\end{equation}

Suppose that the kinetic energy of material points has the form
\begin{equation}
\label{Rvectors}
   \begin{cases}
      T^{(1)}=\frac{m^{(1)}\dot{\vec{x}}^2}{2}=\frac{m^{(1)}}{2}({\partial_a\vec{x}}\cdot{\partial_b\vec{x}})\dot{\xi}^a\dot{\xi}^b;\\
      T^{(2)}=\frac{m^{(2)}\dot{\vec{y}}^2}{2}=\frac{m^{(2)}}{2}({\partial_q\vec{y}}\cdot{\partial_p\vec{y}})\dot{\eta}^q\dot{\eta}^p,
   \end{cases}
\end{equation}
where "$\cdot$" is a scalar product in $E^N$; $a,b,c,d=1\ldots n$ are the indices of the internal variables of the 1-st surfaces;
$q,p,r,s=1\ldots m$ are the indices of the internal variables of the 2-st surfaces.

Potential energy depends only on the distance between the points and is a monotonically increasing function
\begin{equation}
\label{vecr}
\begin{cases}
   U=U(r),\quad r=|\vec{r}|;\\
   \vec{r}=\vec{y}(\eta)-\vec{x}(\xi)=\vec{y}(\eta^1,\cdots,\eta^m)-\vec{x}(\xi^1,\cdots,\xi^n);\\
   |\vec{r}|=\sqrt{\vec{r}^{\ 2}}=\sqrt{r^I r_I}=\sqrt{\sum_{I=1}^N (r^I)^2};\\
   \frac{\partial}{\partial\xi_a}(U(|\vec{r}|))=\frac{\partial U}{\partial r}\frac{\partial r}{\partial\xi_a}=-\frac{1}{r}\frac{\partial U}{\partial r}(\vec{y}-\vec{x})\cdot{\partial_a\vec{x}};\\
   \frac{\partial}{\partial\eta_s}(U(|\vec{r}|))=\frac{\partial U}{\partial r}\frac{\partial r}{\partial\eta_s}=\frac{1}{r}\frac{\partial U}{\partial r}(\vec{y}-\vec{x})\cdot{\partial_s\vec{y}},
\end{cases}
\end{equation}
where $I,J=1\ldots N$ indices enumerate the components of the radius-vectors in Euclidean space.

Then Lagrange function of generalized coordinates ($\xi^1,\cdots,\xi^n,\eta^1,\cdots,\eta^m$) and
generalized velocities ($\dot{\xi}^1,\cdots,\dot{\xi}^n,\dot{\eta}^1,\cdots,\dot{\eta}^m$) has the form
\begin{equation}\label{L0}
   \mathcal{L}_0=T^{(1)}+T^{(2)}-U=
\end{equation}
\[
=\frac{m^{(1)}}{2}({\partial_a\vec{x}}\cdot{\partial_b\vec{x}})\dot{\xi}^a\dot{\xi}^b
+\frac{m^{(2)}}{2}({\partial_q\vec{y}}\cdot{\partial_p\vec{y}})\dot{\eta}^q\dot{\eta}^p
-U(|\vec{y}(\eta)-\vec{x}(\xi)|).
\]

At the surfaces $\mathfrak{N}$, $\mathfrak{M}$ insymbol $g^{(1)},g^{(2)}$ of Riemannian metrics using fundamental quadratic forms of surfaces \cite[(1.1.7),p.3]{MarRat98}
\begin{equation}
\label{Tensors}
\begin{cases}
   g^{(1)}_{ab}=
   m^{(1)}(\partial_{a}\vec{x}\cdot\partial_{b}\vec{x})=m^{(1)}\sum_{I=1}^N \partial_{a} x^I \partial_{b} x_I,\\
   g^{(2)}_{st}=
   m^{(2)}(\partial_{s}\vec{y}\cdot\partial_{t}\vec{y})=m^{(2)}\sum_{J=1}^N \partial_{s} y^J \partial_{t} y_J.
\end{cases}
\end{equation}

\bigskip
\section{ The Lagrangian formalism }
\label{secLagrangeForm}

Taking into account (\ref{Tensors}) Lagrange function (\ref{L0}) has the form:
\begin{equation}\label{L0}
   \mathcal{L}(\xi^1,\ldots,\xi^n,\eta^1,\ldots,\eta^m;\dot{\xi}^1,\ldots,\dot{\xi}^n,\dot{\eta}^1,\ldots,\dot{\eta}^m)=
\end{equation}
\[
=\frac12 g^{(1)}_{ab}\dot{\xi}^a\dot{\xi}^b
+\frac12 g^{(2)}_{st}\dot{\eta}^s\dot{\eta}^t
-U(|\vec{y}(\eta)-\vec{x}(\xi)|).
\]

Let's write down Lagrange equations of our natural mechanical system
\begin{equation}\label{LEq}
\begin{cases}
   \frac{d}{dt}\left(\frac{\partial \mathcal{L}}{\partial\dot{\xi}^a}\right)
  -\frac{\partial \mathcal{L}}{\partial \xi^a}=0;\\
   \frac{d}{dt}\left(\frac{\partial \mathcal{L}}{\partial\dot{\eta}^q}\right)
  -\frac{\partial \mathcal{L}}{\partial \eta^q}=0.
\end{cases}
\end{equation}

For the first equation of (\ref{LEq}) we have
\begin{equation}\label{LEq_Comp}
\begin{cases}
   \frac{\partial\mathcal{L}}{\partial\xi_a}=m^{(1)}({\partial_b\vec{x}}\cdot\partial^2_{ac}\vec{x})\dot{\xi}^c\dot{\xi}^b
      +\partial_a(U(|\vec{y}(\eta)-\vec{x}(\xi)|));\\
   \frac{\partial \mathcal{L}}{\partial\dot{\xi_a}}=m^{(1)}({\partial_a\vec{x}}\cdot{\partial_b\vec{x}})\dot{\xi}^b;\\
   \frac{d}{dt}\left(\frac{\partial \mathcal{L}}{\partial\dot{\xi_a}}\right)=m^{(1)}(({\partial_b\vec{x}}\cdot{\partial^2_{ac}\vec{x}}+{\partial_a\vec{x}}\cdot{\partial^2_{bc}\vec{x}})\dot{\xi^c}\dot{\xi^b}
+({\partial_a\vec{x}}\cdot{\partial_b\vec{x}})\ddot{\xi}^b).
\end{cases}
\end{equation}

Similarly, for the second equation in (\ref{LEq}).

Substituding (\ref{LEq_Comp}) in (\ref{LEq}) we get
\begin{equation}\label{LagrangeEq7}
   \begin{cases}
      g^{(1)}_{ab}\ddot{\xi}^b+
      m^{(1)}({\partial_a\vec{x}}\cdot{\partial^2_{bc}\vec{x}})\dot{\xi^c}\dot{\xi}^b
      +\partial_a(U(|\vec{r}|))=0; \\
      g^{(2)}_{st}\ddot{\eta}^t+
      m^{(2)}({\partial_s\vec{x}}\cdot{\partial^2_{tu}\vec{x}})\dot{\eta}^u\dot{\eta}^t
      +\partial_s(U(|\vec{r}|))=0.
   \end{cases}
\end{equation}

\section{ Lagrangian dynamics on Riemannian manifold }
\label{RimanGeom}

Let’s show that our problem reduces to the dynamics of a point on Riemannian manifold that is the direct product of two Riemannian manifolds (the original surfaces in Euclidean space).

Insymbol
\begin{equation}\label{coord_impals}
   \begin{cases}
      q=(\xi^1,\ldots,\xi^n,\eta^1,\ldots,\eta^m); \\
      \dot{q}=(\dot{\xi}^1,\ldots,\dot{\xi}^n,\dot{\eta}^1,\ldots,\dot{\eta}^m); \\
      p=(\mu_1,\ldots,\mu_n,\nu_1,\ldots,\nu_m).
   \end{cases}
\end{equation}

\begin{equation}\label{metr_matrix}
g=\begin{bmatrix}
g^{(1)}\circ pr_1 & 0 \\
0 & g^{(2)}\circ pr_2
\end{bmatrix},
\end{equation}

\begin{equation}\label{projection}
   \begin{cases}
      pr_1: q=(\xi^1,\ldots,\xi^n,\eta^1,\ldots,\eta^m)\longrightarrow\xi=(\xi^1,\ldots,\xi^n); \\
      pr_2: q=(\xi^1,\ldots,\xi^n,\eta^1,\ldots,\eta^m)\longrightarrow\eta=(\eta^1,\ldots,\eta^m).
   \end{cases}
\end{equation}

Using (\ref{coord_impals}),(\ref{metr_matrix}) lets write (\ref{L0}) in the form
\begin{equation*}\label{L0_metr}
   \mathcal{L}(q,\dot{q})=\frac12 g_{ik}\dot{q}^i\dot{q}^k-U(|\vec{r}(q)|),
\end{equation*}
or
\begin{equation}\label{L0_Marsden}
   \mathcal{L}(q,v)=\frac12\langle v,v\rangle-U(|\vec{r}(q)|),
   \quad v^i=\dot{q}^i, \quad i,j,k,l=1,\dots,n+m.
\end{equation}

Since our configuration space is Riemannian space with metric (\ref{metr_matrix})
that is given by the kinetic energy of the system and Lagrangian has the form (\ref{L0_Marsden})
(see \cite[(7.7.2)--(7.7.3), p.198]{MarRat98}) then our problem can be reduced to the theory of natural mechanical system
with kinetic energy that is determined by the metric of Riemannian manifold \cite[\S 7.7]{MarRat98}.

Let’s transform the equations (\ref{LagrangeEq7}) to the form
\begin{equation}\label{LagrangeEq8}
   \begin{cases}
      g^{(1)}_{ab}\ddot{\xi}^b+
      m^{(1)}({\partial_a\vec{x}}\cdot{\partial^2_{bc}\vec{x}})\dot{\xi^c}\dot{\xi}^b
      =-\partial_a(U(|\vec{r}|)); \\
      g^{(2)}_{st}\ddot{\eta}^t+
      m^{(2)}({\partial_s\vec{x}}\cdot{\partial^2_{tu}\vec{x}})\dot{\eta}^u\dot{\eta}^t
      =-\partial_s(U(|\vec{r}|)).
   \end{cases}
\end{equation}

Then the left-hand sides can be expressed in terms of Riemannian geometry (metric, connection coefficient).

Indeed
\begin{equation}\label{LagrangeInvForm}
g\ddot{q}=\begin{pmatrix}
g^{(1)}\circ pr_1 & 0 \\
\\
\\
0 & g^{(2)}\circ pr_2
\end{pmatrix}
\begin{pmatrix}
\ddot{\xi}^{(1)}\\
\vdots\\
\ddot{\xi}^{(n)}\\
\ddot{\eta}^{(1)}\\
\vdots\\
\ddot{\eta}^{(m)}
\end{pmatrix}
\end{equation}
and
\begin{equation}\label{secondin7}
\begin{split}
   m^{(1)}({\partial_a\vec{x}}\cdot\partial^2_{bc}\vec{x})
&=\frac12(\partial_b g_{ac}+\partial_c g_{ab})-\\
&-\frac12(\partial^2_{ab}\vec{x}\cdot\partial_c\vec{x}+\partial^2_{ac}\vec{x}\cdot\partial_b\vec{x})=\\
&=\frac12(\partial_b g_{ac}+\partial_c g_{ab})-\frac12\partial_a g_{bc}=\\
&=\frac12(\partial_b g_{ac}+\partial_c g_{ab}-\partial_a g_{bc})=\Gamma_{b,ac};
\end{split}
\end{equation}
so far as
\begin{equation*}\label{LagrangeEq5}
   \begin{cases}
      \partial_a g_{bc}=\partial^2_{ab}\vec{x}\cdot\partial_c\vec{x}+\partial_b\vec{x}\cdot\partial^2_{ac}\vec{x};\\
      \partial_b g_{ac}=\partial^2_{ab}\vec{x}\cdot\partial_c\vec{x}+\partial_a\vec{x}\cdot\partial^2_{bc}\vec{x};\\
      \partial_c g_{ab}=\partial^2_{ac}\vec{x}\cdot\partial_b\vec{x}+\partial_a\vec{x}\cdot\partial^2_{bc}\vec{x}.
   \end{cases}
\end{equation*}

Then the equations (\ref{LagrangeEq8}) take the form
\begin{equation}\label{LagrangeEq9}
   \begin{cases}
      g^{(1)}_{ab}\ddot{\xi}^b+
      \Gamma_{b,ac}\dot{\xi^c}\dot{\xi}^b
      =-\partial_a(U(|\vec{r}|)); \\
      g^{(2)}_{st}\ddot{\eta}^t+
      \Gamma_{t,su}\dot{\eta}^t\dot{\eta}^u
      =-\partial_s(U(|\vec{r}|));\\
      \Gamma_{b,ac}=\frac12(\partial_b g_{ac}+\partial_c g_{ab}-\partial_a g_{bc});\\
      \Gamma_{t,su}=\frac12(\partial_t g_{su}+\partial_u g_{st}-\partial_s g_{tu}).
   \end{cases}
\end{equation}

Let's use the well-known relation for the Christoffel symbol
\begin{equation}\label{SymCryst1}
   \Gamma^a_{bc}=g^{ad}\Gamma_{b,dc},
\end{equation}
where $g^{ad}$ is contravariant metric tensor given by an inverse matrix with respect to a covariant matrix. I.e.
\begin{equation}\label{MetricRev}
   \begin{cases}
      g_{(1)}^{ab}g^{(1)}_{bc}=\delta^a_c;\\
      g_{(2)}^{st}g^{(2)}_{tu}=\delta^s_u;\\
      g_{ik}g^{kl}=\delta^l_i,
   \end{cases}
\end{equation}
where $g$ in last relation is taken from (\ref{metr_matrix}) with taking into account the block structure of this matrix.

Then the Lagrange equations of motion take the form
\begin{equation}\label{LagrangeEq6}
   \begin{cases}
      \ddot{\xi}^a+(\Gamma^{(1)})^a_{bc}\dot{\xi}^b\dot{\xi}^c=-g_{(1)}^{ac}\partial_c(U(|\vec{r}|));\\
      \ddot{\eta}^s+(\Gamma^{(2)})^s_{tu}\dot{\eta}^t\dot{\eta}^u=-g_{(2)}^{st}\partial_t(U(|\vec{r}|)),
   \end{cases}
\end{equation}
i.e.
\begin{equation}\label{LagrangeSymCryst3}
   \ddot{q}^i+\Gamma^i_{jk}\dot{q}^j\dot{q}^k+g^{ik}\partial_k(U(|\vec{r}|))=0.
\end{equation}

\bigskip
\section{ The Hamiltonian equations and dissipative term }
\label{secDalamber}

Let's represent the equation (\ref{LagrangeSymCryst3}) in the form
\begin{equation}\label{LagrangeEqCryst4}
   \ddot{q}^i=\gamma(q,\dot{q})^i-(\nabla U)^i,
\end{equation}
where
\begin{equation}\label{Grad}
\begin{cases}
   \gamma(q,\dot{q})^i=-\Gamma^i_{jk}\dot{q}^j\dot{q}^k;\\
   (\nabla U)^i=g^{ik}\frac{\partial U}{\partial q^k}.
\end{cases}
\end{equation}

Equation (\ref{LagrangeEqCryst4}) fully coincides with \cite[(7.7.3),p.198]{MarRat98}.

The transition to the first-order Hamiltonian equations is realized with the help of Legendre transformations
\begin{equation}\label{LegandreTransf1}
\begin{cases}
   p_i=v_i=g_{ij}v^j=g_{ij}\dot{q}^j;\\
   \mathcal{H}(q,p)=E(v)=A(v)-\mathcal{L}(v)=\\
   =\langle v,v\rangle_q -\mathcal{L}(v)=\frac12\langle v,v\rangle_q+U(q).
\end{cases}
\end{equation}
\begin{equation}\label{LegandreTransf2}
\begin{cases}
   \dot{q}=v;\\
   \dot{v}=\gamma(q,v)-\nabla U(q).
\end{cases}
\end{equation}
\begin{remark}
Strictly speaking the Hamiltonian equations must be written with respect to the covariant vector $p_i$
rather than a contravariant vector $v^i$. However, the equations of the first order (\ref{LegandreTransf2})
are completely equivalent to the Hamiltonian equations that can easily be obtained from (\ref{LegandreTransf2})
by using covariant differentiation.
\end{remark}

Indeed
\begin{equation}\label{CovarDeriv0}
\begin{cases}
      \frac{d v}{dt}-\gamma(q,v)=\frac{D v}{dt};\\
      \frac{D v}{dt} = -\nabla U;\\
      g^\flat\frac{D v}{dt} = -d U;\\
            \frac{D g}{dt} = 0;\\
      g^\flat\frac{D v}{dt} = \frac{D (g^\flat v)}{dt} = \frac{D p}{dt};\\
      \frac{D p}{dt} = - d U,
   \end{cases}
\end{equation}
where $\frac{D}{dt}$ is a symbol of the covariant derivative along the trajectory.

As it shown in \cite[с.205]{MarRat98} the dissipation can be described via the Rayleigh function and the generalized Lagrange equation takes the form
\begin{equation}\label{Rayleigh1}
   \frac{d}{dt}\left(\frac{\partial \mathcal{L}}{\partial \dot{q}^i}\right)-\frac{\partial \mathcal{L}}{\partial q^i}=-\frac{\partial \mathcal{R}}{\partial \dot{q}^i},
\end{equation}
where
\[
\mathcal{R}=\frac12 R_{ij}(q)v^i v^j.
\]
Then the equations (\ref{LegandreTransf2}) taking into account the dissipation can be written
\begin{equation}\label{Rayleigh2}
\begin{cases}
   \dot{q}=v;\\
   \dot{v}=\gamma(q,v)-\nabla U(q)-F_\mathcal{R}(q,v);\\
    F_\mathcal{R}^i=g^{ik}\frac{\partial \mathcal{R}}{\partial v^k}
\end{cases}
\end{equation}
\begin{remark}
First the Lagrange equations (\ref{LagrangeEq7}) for our system were written in the local coordinate system.
However, the principle of action for the Lagrange equation (\ref{LegandreTransf2}) and d'Alembert principle \cite[p.202--203]{MarRat98} from that follows the equations (\ref{Rayleigh2}) that formulated on Riemannian manifold,
independently of local coordinate systems. Thus, the equations (\ref{Rayleigh2})
define a global vector field on the Riemannian manifold of our natural mechanical system.
\end{remark}

\bigskip
\section{ Liapounov function }
\label{Lyapunov}

\begin{lemma}
The point $(q_0,0)$ where $\nabla U(q_0)=0$ and $v=0$ is an equilibrium point for the (\ref{Rayleigh2}).
\end{lemma}

This is obvious from the form of the right-hand sides of the equations (\ref{Rayleigh2}).

\begin{proposition}
The function $L(q,v)=E(q,v)-U(q_0)$ is the Lyapunov function for the field (\ref{Rayleigh2}).
\end{proposition}

\begin{itemize}
\item $L(q_0,0)=0$ by construction.
\item For provide $L(q,v)>0$ for anyone $(q,v)\in U\backslash \{(q_0,0)\}$ it is necessary that $U(q)>U(q_0)$ on a whole neighborhood area (as known that the kinetic energy is positive definite).
\item Let’s show that $\dot{L}(q,v):=X[L](q,v)<0$ anywhere $(q,v)\in U\backslash \{(q_0,0)\}$.
\end{itemize}
\begin{equation*}\label{Energy1}
\begin{split}
   \frac{dE}{dt}=\frac{DE}{dt}=\frac{D}{dt}\left(\frac12\langle v,v\rangle_q+U(q)\right)
&=\left\langle v,\frac{Dv}{dt}\right\rangle+\frac{dU}{dt}=\\
&=-\langle v,\nabla U\rangle-\langle v,\nabla^{(v)}\mathcal{R}\rangle+\langle\nabla U, v\rangle=\\
&=-\langle v,\nabla^{(v)}\mathcal{R}\rangle=\\
&=-v^i\frac{\partial \mathcal{R}}{\partial v^i}=-2\mathcal{R}<0,
\end{split}
\end{equation*}
where
$E(q,v)$ is the total energy of the Lagrangian system without dissipation.

In last relation we use the fact that Rayleigh function is a quadratic positive-definite form of  velocities,
i.e. homogeneous velocity function of 2-degree.

Then from Lyapunov theorem (the second method)\cite{OrtegaLect,ModContr} take place the asymptotic stability of the equations (\ref{Rayleigh2}) at the equilibrium point.

\begin{remark}
Let {\it note} that the potential energy of type $U=k r^2$, where $k>0$ is suitable case.
\end{remark}

\end{document}